\documentclass[11pt]{article}

\usepackage{amssymb, amsmath, amsthm, latexsym}
\usepackage{fullpage, color}
\usepackage{url}
\usepackage{algorithm}
\usepackage{algpseudocode}
\usepackage{tabularx}
\usepackage{paralist}
\usepackage{mathtools}
\usepackage{tcolorbox}
\usepackage{xcolor}

\usepackage{bbm} 

\usepackage{makecell}
\usepackage{multirow}
\usepackage{booktabs}

\usepackage[flushleft]{threeparttable} 

\usepackage{caption}
\usepackage{multirow}
\usepackage{colortbl}
\definecolor{bgcolor}{rgb}{0.8,1,1}
\definecolor{bgcolor2}{rgb}{0.8,1,0.8}
\definecolor{bgcolor3}{rgb}{0.67, 0.94, 0.82}

\usepackage{nicefrac}       

\definecolor{niceblue}{rgb}{0.0,0.19,0.56}

\usepackage{hyperref}
\hypersetup{colorlinks,linkcolor={blue},citecolor={niceblue},urlcolor={blue}}

\usepackage[numbers]{natbib}
\usepackage{sidecap}

\usepackage{pifont}
\usepackage{subfigure}

\newcommand{\R}{\mathbb{R}}

\def\<#1,#2>{\left\langle #1,#2\right\rangle}

\usepackage{mdframed} 
\usepackage{thmtools}

\definecolor{shadecolor}{gray}{0.9}
\declaretheoremstyle[
headfont=\normalfont\bfseries,
notefont=\mdseries, notebraces={(}{)},
bodyfont=\normalfont,
postheadspace=0.5em,
spaceabove=1pt,
mdframed={
  skipabove=8pt,
  skipbelow=8pt,
  hidealllines=true,
  backgroundcolor={shadecolor},
  innerleftmargin=4pt,
  innerrightmargin=4pt}
]{shaded}

\declaretheorem[style=shaded,within=section]{definition}
\declaretheorem[style=shaded,sibling=definition]{theorem}
\declaretheorem[style=shaded,sibling=definition]{proposition}
\declaretheorem[style=shaded,sibling=definition]{assumption}
\declaretheorem[style=shaded,sibling=definition]{corollary}

\declaretheorem[style=shaded,sibling=definition]{lemma}

\usepackage[colorinlistoftodos,bordercolor=orange,backgroundcolor=orange!20,linecolor=orange,textsize=scriptsize]{todonotes}

\newcommand{\cC}{{\cal C}}
\newcommand{\cD}{{\cal D}}

\newcommand{\cL}{{\cal L}}

\newcommand{\cO}{{\cal O}}

\newcommand{\mA}{{\bf A}}

\newcommand{\EE}{\mathbf{E}}

\def\R{\mathbb{R}}

\def\R{\mathbb R}

\def\EE{\mathbb E}
\def\PP{\mathbb P}

\def\tx{\tilde{x}}

\usepackage{xspace}

\newcommand{\algname}[1]{{\sf #1}\xspace}

\usepackage{hyperref}

\usepackage{accents}
\newlength{\dhatheight}

\begin{document}
\title{\textbf{Distributed Methods with Absolute Compression}\\ \textbf{and Error Compensation}\thanks{The research was supported by Russian Science Foundation grant (project No. 21-71-30005).}}
%
%
\author{Marina Danilova$^{1,2}$ \and
Eduard Gorbunov$^2$}
\date{$^1$Institute of Control Sciences of RAS, Moscow, Russia\\ $^2$Moscow Institute of Physics and Technology, Moscow, Russia}
\maketitle              
\begin{abstract}
Distributed optimization methods are often applied to solving huge-scale problems like training neural networks with millions and even billions of parameters. In such applications, communicating full vectors, e.g., (stochastic) gradients, iterates, is prohibitively expensive, especially when the number of workers/nodes is large. Communication compression is a powerful approach to alleviating this issue, and, in particular, methods with biased compression and \textit{error compensation} are extremely popular due to their practical efficiency. Sahu et al. (2021) \cite{sahu2021rethinking} propose a new analysis of Error Compensated \algname{SGD} (\algname{EC-SGD}) for the class of absolute compression operators showing that in a certain sense, this class contains optimal compressors for \algname{EC-SGD}. However, the analysis was conducted only under the so-called $(M,\sigma^2)$-bounded noise assumption. In this paper, we generalize the analysis of \algname{EC-SGD} with absolute compression to the arbitrary sampling strategy and propose the first analysis of Error Compensated Loopless Stochastic Variance Reduced Gradient method (\algname{EC-LSVRG}) \cite{gorbunov2020linearly}  with absolute compression for (strongly) convex problems. Our rates improve upon the previously known ones in this setting. Numerical experiments corroborate our theoretical findings.
\end{abstract}

\section{Introduction}
In the recent few years, distributed optimization methods has been receiving a lot of attention from various research communities and, especially, from the machine learning one. This can be explained by the need of training deep learning models with billions of parameters on the hundreds of gigabytes of data \cite{brown2020language} (and sometimes even this is not enough \cite{kaplan2020scaling}). Clearly, such problems cannot be solved in a reasonable time on a single yet powerful machine \cite{li2020openai}. Next, distributed methods are literally the only possible choices in such applications like Federated Learning (FL) \cite{konecny2016federated,FEDLEARN,kairouz2019advances}, where the data is privately stored on multiple devices.

Due to the huge dimensions of corresponding problems and large number of workers in the networks na\"ive methods like centralized Parallel \algname{SGD} \cite{zinkevich2010parallelized} suffer from the so-called \emph{communication bottleneck}. This phenomenon means that a method spends much more time on communication rounds than on computations. A natural and popular way of addressing this issue is \emph{communication compression} \cite{seide20141} -- a technique that uses special compression operators called \emph{compressors} applied to the information that devices send through the network. 

The works on distributed methods with compression usually focus either on unbiased compressors \cite{alistarh2017qsgd,mishchenko2019distributed,horvath2019stochastic} like RandK or $\ell_2$-quantization (e.g., see \cite{beznosikov2020biased}) or on biased compressors \cite{seide20141,stich2018sparsified,beznosikov2020biased,gorbunov2020linearly} like TopK. Although the world of unbiased compressors has richer theory, the methods with biased compressors are very popular due to their efficiency in practice. However, to make them convergent one has to apply additional tricks on top of \algname{SGD}, e.g., \emph{error-compensation} \cite{seide20141,stich2018sparsified,stich2020error,beznosikov2020biased}.

Error Compensated \algname{SGD} (\algname{EC-SGD}) was proposed in \cite{seide20141} where the authors demonstrated its efficiency in practical tasks, but the first theoretical analysis of \algname{EC-SGD} was given in \cite{stich2018sparsified} and tightened in \cite{stich2020error}. This analysis was extended in various directions including (but not limited to) decentralized communications \cite{koloskova2019decentralized}, arbitrary sampling and variance reduction (with the first linearly convergent variants) \cite{gorbunov2020linearly}, acceleration \cite{qian2020error}, and also some prominent alternatives were proposed \cite{horvath2019stochastic,richtarik2021ef21}. However, the compression operators in these papers are usually assumed to be \emph{$\delta$-contractive}.\footnote{The mapping (possibly stochastic) $\cC : \R^d \to \R^d$ is called \emph{$\delta$-contractive} compressor if there exists $\delta \in (0,1]$ such that $\EE[\|\cC(x) - x\|^2] \leq (1-\delta)\|x\|^2$ for all $x \in\R^d$.}

Recently, the authors of \cite{sahu2021rethinking} developed a new analysis of \algname{EC-SGD} with \emph{absolute compressors} \cite{tang2019doublesqueeze,sahu2021rethinking}, i.e., such (stochastic) operators $\cC$ that for some $\Delta \geq 0$ the inequality $\EE[\|\cC(x) - x\|^2] \leq \Delta^2$ holds for all $x \in \R^d$, where $\EE[\cdot]$ denotes an expectation.. In particular, they proved that this class contains special operators called hard-threshold sparsifiers that are optimal in view of total error minimization (a special quantity arising in the analysis of \algname{EC-SGD}) for any \emph{fixed} sequence of errors. Moreover, the authors of \cite{sahu2021rethinking} derived convergence rates for \algname{EC-SGD} with absolute compressors under $(M,\sigma^2)$-bounded noise assumption and illustrate the theoretical and practical benefits of absolute compressors compared to $\delta$-contractive ones. However, several fruitful directions were unexplored for \algname{EC-SGD} with absolute compressor including more general analysis of the standard version of the method and variants with variance reduction.

\subsection{Main Contributions}
%
\begin{itemize}
    \item[$\diamond$] \textbf{Unified analysis of \algname{EC-SGD} with absolute compressors.} We propose a generalized analysis of \algname{EC-SGD} with absolute compression covering different stochastic estimators under various assumptions. In particular, we consider the simplified version of the parametric assumption from \cite{gorbunov2020linearly} (see Assumption~\ref{as:key_assumption} and the discussion after) and derive a general result on the convergence of \algname{EC-SGD} with absolute compressors (Theorem~\ref{thm:main_result}). The considered assumption covers various setups including the one from \cite{sahu2021rethinking} and the derived result gives sharp rates.
     
    \item[$\diamond$] \textbf{\algname{EC-SGD} with absolute compression and arbitrary sampling.} To illustrate the flexibility of our approach, we propose the first analysis of \algname{EC-SGD} with absolute compression and arbitrary sampling. The derived bounds are superior to the ones from \cite{sahu2021rethinking} under certain assumptions.
    
    \item[$\diamond$] \textbf{\algname{EC-LSVRG} with absolute compression and arbitrary sampling.} As a special case of our theoretical framework, we obtain the analysis of \algname{EC-LSVRG} \cite{gorbunov2020linearly} with absolute compression. Moreover, in contrast to \cite{gorbunov2020linearly}, we handle non-uniform sampling in \algname{EC-LSVRG}. The derived rate for \algname{EC-LSVRG} with absolute compression has the leading term proportional to $\nicefrac{1}{K^2}$, while the leading term for \algname{EC-SGD} is proportional to $\nicefrac{1}{K}$, where $K$ is the total number of iterations. 
    
    \item[$\diamond$] \textbf{Numerical experiments.} We conduct several numerical experiments to support our theory and compare the performance of \algname{EC-LSVRG} with hard-threshold and TopK sparsifiers. The numerical results corroborate our theoretical findings and highlight the benefits of using absolute compressors.
\end{itemize}

\subsection{Preliminaries}
\begin{algorithm}[t]
\caption{Error-Compensated Stochastic Gradient Descent (\algname{EC-SGD})}
\label{alg:EC_SGD}   
\begin{algorithmic}[1]
\Require starting point $x^0$, stepsize $\gamma > 0$, number of iterations $K \ge 0$
\State Set $e_i^0 = 0$ for all $i = 1,\ldots, n$
\For{$k=0,\ldots, K-1$}
\State Server broadcasts $x^k$ to all workers
\For{$i = 1,\ldots, n$ in parallel \textbf{do}}
\State Compute stochastic gradient $g_i^k$ and send $v_i^k = \gamma \cC\left(\tfrac{e_i^k + \gamma g_i^k}{\gamma}\right)$ to the server
\State Update error-vector: $e_i^{k+1} = e_i^k + \gamma g_i^k - v_i^k$
\EndFor 
\State Server gathers $v_1^k, v_2^k, \ldots, v_n^k$ from all workers and computes $v^k = \tfrac{1}{n}\sum_{i=1}^n v_i^k$
\State Set $x^{k+1} = x^k - v^k$
\EndFor 
\end{algorithmic}
\end{algorithm}

\noindent\textbf{Problem.} We consider a classical centralized optimization problem
\begin{equation}
    \min\limits_{x\in \R^d}\left\{f(x) = \tfrac{1}{n}\sum_{i=1}^n f_i(x)\right\}, \label{eq:distributed_problem}
\end{equation}
where the information defining differentiable functions $f_1,\ldots, f_n:\R^d \to \R$ is distributed among $n$ workers/clients/devices connected with parameter-server in a centralized way. In particular, (stochastic) gradients of function $f_i$ are available to client $i$ only. Throughout the work, we assume that the solution $x^*$ of problem \eqref{eq:distributed_problem} is unique and that the function $f$ is convex and $\mu$-quasi strongly convex \cite{necoara2019linear}, where the later is a relaxation of strong convexity meaning that
\begin{equation}
    \forall x \in \R^d\quad f(x^*) \ge f(x) + \langle\nabla f(x), x^* - x \rangle + \tfrac{\mu}{2}\|x - x^*\|^2,\quad \mu \geq 0. \label{eq:quasi_strong_cvx}
\end{equation}
Moreover, we assume that $f(x)$ is $L$-smooth, i.e., $\|\nabla f(x) - \nabla f(y)\| \leq L\|x-y\|$ holds for all $x,y \in \R^d$.

\noindent\textbf{Compression.} In this work, we focus on \emph{absolute compression operators} \cite{tang2019doublesqueeze,sahu2021rethinking}.
\begin{definition}\label{def:absolute_compression}
    The mapping (possibly stochastic) $\cC : \R^d \to \R^d$ is called absolute compression operator/absolute compressor if there exists $\Delta \geq 0$ such that
    \begin{equation}
        \EE\left[\|\cC(x) - x\|^2\right] \leq \Delta^2,\quad \forall x \in \R^d. \label{eq:absolute_compression_def}
    \end{equation}
\end{definition}

\noindent An example of absolute compressor is hard-threshold sparsifier $\cC_{\text{HT}}(x)$ \cite{strom2015scalable,dutta2020discrepancy,sahu2021rethinking} defined as $[\cC_{\text{HT}}(x)]_i = [x]_i$ if $|[x]_i| \geq \lambda$ and $[\cC_{\text{HT}}(x)]_i = 0$ otherwise  
 for some $\lambda \geq 0$, where $[\cdot]_i$ denotes the $i$-th component of the vector. One can show that $\cC_{\text{HT}}(x)$ satisfies \eqref{eq:absolute_compression_def} with $\Delta = \lambda\sqrt{d}$. Other examples include (stochatsic) rounding schemes with bounded error \cite{gupta2015deep} and scaled integer rounding \cite{sapio2021scaling}.

\noindent \textbf{Paper organization.} In Section~\ref{sec:unif_analysis}, we formulate our general result on the convergence of \algname{EC-SGD} with absolute compression. Next, we provide particular examples of the variants of \algname{EC-SGD} fitting our framework -- \algname{EC-SGD} with arbitrary sampling (Section~\ref{sec:arb_sampling}) and \algname{EC-LSVRG} (Section~\ref{sec:lsvrg}) -- and discuss the convergence guarantees obtained for them. Finally, in Section~\ref{sec:num_exp}, we discuss the results of numerical experiments supporting our theoretical findings. The proofs are delegated to the Appendix.

\section{Unified Analysis}\label{sec:unif_analysis}

In our analysis, we rely on a simplified version of Assumption~3.3 from\footnote{The unified analysis of stochastic first-order methods was proposed in \cite{gorbunov2019unified} for quasi strongly convex problems. After that, this idea was extended to the case of convex functions \cite{khaled2020unified}, methods with error feedback \cite{gorbunov2020linearly} and local updates \cite{gorbunov2021local}, and to the methods for solving variational inequalities and min-max problems \cite{gorbunov2021stochastic,beznosikov2022stochastic}.} \cite{gorbunov2020linearly}.

\begin{assumption}[Key Parametric Assumption]\label{as:key_assumption}
    For all $k \ge 0$ the average of stochastic gradients used in \algname{EC-SGD} is an unbiased estimate of $\nabla f(x^k)$, i.e., for $g^k = \tfrac{1}{n}\sum_{i=1}^n g_i^k$ we have $\EE_k[g^k] = \nabla f(x^k)$ for all $k \ge 0$, where $\EE_k[\cdot]$ denotes an expectation w.r.t.\ the randomness coming from iteration $k$. Moreover, there exist non-negative parameters $A, B, C, D_1, D_2 \geq 0$, $\rho \in (0,1]$, and sequence of (possibly random) variables $\{\sigma_k^2\}_{k\ge 0}$ such that for all $k \ge 0$ the iterates produced by \algname{EC-SGD} and the objective function $f$ satisfy
    \begin{eqnarray}
        \EE_k\left[\|g^k\|^2\right] &\leq& 2A\left(f(x^k) - f(x^*)\right) + B\sigma_k^2 + D_1, \label{eq:second_moment_bound_general}\\
        \EE_k\left[\sigma_{k+1}^2\right] &\leq& (1-\rho)\sigma_k^2 + 2C\left(f(x^k) - f(x^*)\right) + D_2. \label{eq:sigma_k_bound_geenral}
    \end{eqnarray}
\end{assumption}

As it is shown in \cite{gorbunov2020linearly}, the above assumption is very general and covers various algorithms in different settings. In Section~\ref{sec:arb_sampling} and \ref{sec:lsvrg}, we consider two particular examples when Assumption~\ref{as:key_assumption} is satisfied. In all known special cases, parameters $A$ and $C$ are typically related to the smoothness properties of the problem, $\sigma_k^2$ describes the variance reduction process (with ``rate'' $\rho$), $D_1$ and $D_2$ are remaining noises not handled by variance reduction, and $B$ is some constant.

Under Assumption~\ref{as:key_assumption} we derive the following result in Appendix~\ref{appendix:main_proofs}.
\begin{theorem}\label{thm:main_result}
    Let function $f$ be convex, $\mu$-quasi strongly convex (with unique solution $x^*$), $L$-smooth, and Assumption~\ref{as:key_assumption} hold. Assume that $0 < \gamma \leq \nicefrac{1}{4(A + CF)}$, where $F = \nicefrac{4B}{3\rho}$. Then, for all $K \ge 0$ the iterates produced by \algname{EC-SGD} with absolute compression operator $\cC$ (see Definition~\ref{def:absolute_compression}) satisfy
    \begin{eqnarray}
        \EE\left[f(\overline{x}^K) - f(x^*)\right] &\leq&  \tfrac{(1\!-\!\eta)^{K\!+\!1} 2\EE[T_0]}{\gamma}\! +\! 2\gamma\!\left(D_1\! +\! FD_2\! + 3L\gamma \Delta^2\right)\!, \text{ if } \mu > 0, \label{eq:main_result_str_cvx}\\
        \EE\left[f(\overline{x}^K) - f(x^*)\right] &\leq& \tfrac{2\EE[T_0]}{\gamma (K+1)}\! +\! 2\gamma\!\left(D_1\! + FD_2\! + 3L\gamma \Delta^2\right)\!, \text{ if } \mu = 0, \label{eq:main_result_cvx}
    \end{eqnarray}
    where $\overline{x}^K = \tfrac{1}{W_K}\sum_{k=0}^{K} w_k x^k$, $w_k = (1 - \eta)^{-(k+1)}$, $\eta = \min\{\nicefrac{\gamma\mu}{2}, \nicefrac{\rho}{4}\}$, $W_K = \sum_{k=0}^{K} w_k$, and $T_0 = \|x^0 - x^*\|^2 + F\gamma^2\sigma_0^2$
\end{theorem}
Upper bounds \eqref{eq:main_result_str_cvx} and \eqref{eq:main_result_cvx} establish convergence to some neighborhood of the solution (in terms of the functional suboptimality). Applying Lemmas~\ref{lem:str_cvx_complexity} and \ref{lem:cvx_complexity}, we derive the convergence rates to the exact optimum.

\begin{corollary}\label{cor:main_result}
    Let the assumptions of Theorem~\ref{thm:main_result} hold. Then, there exist choices of the stepsize $\gamma$ 
     such that \algname{EC-SGD} with absolute compressor guarantees $\EE[f(\overline{x}^K) - f(x^*)]$ of the order
    \begin{eqnarray}
        \widetilde{\cO}\left(\!(A\!+\!CF)\EE[\hat T_0]\exp\left(\!-\min\left\{\tfrac{\mu}{A\!+\!CF}, \rho\right\}\!K\!\right)\! +\! \tfrac{D_1\! +\! FD_2}{\mu K}\! +\! \tfrac{L\Delta^2}{\mu^2 K^2}\!\right)\!,\; (\mu > 0) \label{eq:complexity_str_cvx}\\
        \cO\left(\tfrac{(A\!+\!CF)R_0^2}{K} + \tfrac{R_0^2 \sqrt{B\EE[\sigma_0^2]}}{K\sqrt{\rho}} + \sqrt{\tfrac{R_0^2(D_1 + FD_2)}{K}} + \tfrac{L^{\nicefrac{1}{3}}R_0^{\nicefrac{4}{3}}\Delta^{\nicefrac{2}{3}}}{K^{\nicefrac{2}{3}}}\right)\!,\; (\mu = 0) \label{eq:complexity_cvx}
    \end{eqnarray}
    when $\mu = 0$, where $R_0 = \|x^0 - x^*\|$, $\hat T_0 = R_0^2 + \tfrac{F}{16(A+CF)^2}\sigma_0^2$.
\end{corollary}

This general result allows to obtain the convergence rates for all methods satisfying Assumption~\ref{as:key_assumption} via simple plugging of parameters $A, B, C, D_1, D_2$ and $\rho$ in the upper bounds \eqref{eq:complexity_str_cvx} and \eqref{eq:complexity_cvx}. For example, due to such a flexibility we recover the results from \cite{sahu2021rethinking}, where the authors assume that each $f_i$ is convex and $L_i$-smooth, $f$ is $\mu$-quasi strongly convex\footnote{Although the authors of \cite{sahu2021rethinking} write in Section 3.1 that all $f_i$ are $\mu$-strongly convex, in the proofs, they use convexity of $f_i$ and quasi-strong monotonicity of $f$.}, and stochastic gradients have $(M,\sigma^2)$-bounded noise, i.e., $\EE_k[\|g_i^k - \nabla f_i(x^k)\|^2] \leq M\|\nabla f_i(x^k)\|^2 + \sigma^2$ for all $i \in [n]$. In the proof of their main results (see inequality (24) in \cite{sahu2021rethinking}), they derive an upper-bound for $\EE_k[\|g^k\|^2]$ implying that Assumption~\ref{as:key_assumption} is satisfied in this case with $A = L + \tfrac{M\max_{i\in [n]}L_i}{n}$, $B = 0$, $\sigma_{k}^2\equiv 0$, $D_1 = \tfrac{2M\zeta_*^2 + \sigma^2}{n}$, $C = 0$, $D_2 = 0$, where $\zeta_*^2 = \tfrac{1}{n}\sum_{i=1}^n \|\nabla f_i(x^*)\|^2$ measures the heterogeneity of local loss functions at the solution. Plugging these parameters in Corollary~\ref{cor:main_result} we recover\footnote{When $\mu = 0$, our result is tighter than the corresponding one from \cite{sahu2021rethinking}.} the rates from \cite{sahu2021rethinking}. In particular, when $\mu > 0$ the rate is
\begin{eqnarray}
    \widetilde{\cO}\left(\!(L+\tfrac{M\max_{i\in [n]}L_i}{n})R_0^2\exp\left(\!-\tfrac{\mu}{L+\tfrac{M\max_{i\in [n]}L_i}{n}}K\!\right) + \tfrac{2M\zeta_*^2 + \sigma^2}{\mu nK} + \tfrac{L\Delta^2}{\mu^2 K^2}\!\right)\!. \label{eq:complexity_str_cvx_EC_SGD_M_sigma}
\end{eqnarray}
Below we consider two other examples when Assumption~\ref{as:key_assumption} is satisfied.

\section{Absolute Compression and Arbitrary Sampling}\label{sec:arb_sampling}
Consider the case when each function $f_i$ has a finite-sum form, i.e.,  $f_i(x) = \tfrac{1}{m}\sum_{j=1}^m f_{ij}(x),$ 
 which is a classical situation in distributed machine learning. Typically, in this case, workers sample (e.g., uniformly with replacement) some batch of functions from their local finite-sums to compute the stochastic gradient. To handle a wide range of sampling strategies, we follow \cite{gower2019sgd} and consider a stochastic reformulation of the problem:
\begin{equation}
    f(x) = \EE_{\xi \sim \cD}\left[f_\xi(x)\right],\quad f_\xi(x) = \tfrac{1}{n}\sum_{i=1}^n f_{\xi_i}(x),\quad f_{\xi_i}(x) = \tfrac{1}{m}\sum_{j=1}^m \xi_{ij} f_{ij}(x), \label{eq:stoch_reformulation}
\end{equation}
where $\xi = (\xi_1^\top,\ldots, \xi_n^\top)$ and random vector $\xi_i = (\xi_{i1},\ldots, \xi_{im})^\top$ defines the sampling strategy with distribution $\cD_i$ such that $\EE[\xi_{ij}] = 1$ for all $i\in [n]$, $j \in [m]$. We assume that functions $f_{\xi_i}(x)$ satisfy \emph{expected smoothness} property \cite{gower2019sgd,gorbunov2020linearly}.
\begin{assumption}[Expected Smoothness]\label{as:expected_smoothness}
    Functions $f_1,\ldots, f_n$ are $\cL$-smooth in expectation w.r.t.\ distributions $\cD_1,\ldots, \cD_n$. That is, there exists constant $\cL > 0$ such that for all $x \in \R^d$ and for all $i = 1,\ldots, n$
    \begin{equation}
        \EE_{\xi_i \sim \cD_i}\left[\left\|\nabla f_{\xi_i}(x) - \nabla f_{\xi_i}(x^*)\right\|^2\right] \leq 2\cL\left(f_i(x) - f_i(x) - \langle \nabla f_i(x^*), x - x^* \rangle\right). \notag
    \end{equation}
\end{assumption}
One can show that this assumption (and reformulation itself) covers for a wide range of situations \cite{gower2019sgd}. For example, when all functions $f_{ij}$ are convex and $L_{ij}$-smooth, then for the classical uniform sampling we have $\PP\{\xi_{i} = m e_{j}\} = \tfrac{1}{m}$ and $\cL = \cL_{\text{US}} = \max_{i\in [n], j\in [m]} L_{ij}$, where $e_j\in \R^m$ denotes the $j$-th vector in the standard basis in $\R^m$. Moreover, importance sampling $\PP\{\xi_{i} = \tfrac{m \overline{L}_i}{L_{ij}} e_{j}\} = \tfrac{L_{ij}}{m\overline{L}_i}$, where $\overline{L}_i = \tfrac{1}{m}\sum_{i=j}^m L_{ij}$, also fits Assumption~\ref{as:expected_smoothness} with $\cL = \cL_{\text{IS}} = \max_{i\in [n]}\overline{L}_i,$ which can be significantly smaller than $\cL_{\text{US}}$.

Next, it is worth mentioning that Assumption~\ref{as:expected_smoothness} and $(M,\sigma^2)$-bounded noise assumption used in \cite{sahu2021rethinking} cannot be compared directly, i.e., in general, none of them is stronger than another. However, in contrast to $(M,\sigma^2)$-bounded noise assumption, Assumption~\ref{as:expected_smoothness} is satisfied whenever $f_{ij}(x)$ are convex and smooth.

Consider a special case of \algname{EC-SGD} applied to the stochastic reformulation \eqref{eq:stoch_reformulation}, i.e., let $g_i^k = \nabla f_{\xi_i^k}(x^k)$, where $\xi_i^k$ is sampled from $\cD_i$ independently from previous steps and other workers. Since this version of \algname{EC-SGD} supports arbitrary sampling we will call it \algname{EC-SGD-AS}. In this setup, we show\footnote{Proposition~\ref{prop:EC_SGD_AS_fits_assumption} is a refined version of Lemma~J.1 from \cite{gorbunov2020linearly}.} that \algname{EC-SGD-AS} fits Assumption~\ref{as:key_assumption} (the proof is deferred to Appendix~\ref{appendix:proofs_special_cases}).
\begin{proposition}\label{prop:EC_SGD_AS_fits_assumption}
    Let $f$ be $L$-smooth, $f_i$ have finite-sum form, and Assumption~\ref{as:expected_smoothness} hold. Then the iterates produced by \algname{EC-SGD-AS} satisfy Assumption~\ref{as:key_assumption} with $A = L + \nicefrac{2\cL}{n}$, $B = 0$, $D_1 = \tfrac{2\sigma_*^2}{n} = \tfrac{2}{n^2}\sum_{i=1}^n \EE[\|\nabla f_{\xi_i}(x^*) - \nabla f_i(x^*)\|^2]$, $\sigma_k^2 \equiv 0$, $\rho = 1$, $C = 0$, $D_2 = 0$.
\end{proposition}

Plugging the parameters from the above proposition in Theorem~\ref{thm:main_result} and Corollary~\ref{cor:main_result}, one can derive convergence guarantees for \algname{EC-SGD-AS} with absolute compression operator. In particular, our general analysis implies the following result.
\begin{theorem}\label{thm:EC_SGD_AS_convergence}
   Let the assumptions of Proposition~\ref{prop:EC_SGD_AS_fits_assumption} hold. Then, there exist choices of stepsize $0 < \gamma \leq (4L + \nicefrac{8\cL}{n})^{-1}$ such that \algname{EC-SGD-AS} with absolute compressor guarantees $\EE[f(\overline{x}^K) - f(x^*)]$ of the order
    \begin{eqnarray}
        \widetilde{\cO}\left(\left(L + \tfrac{\cL}{n}\right)R_0^2\exp\left(-\tfrac{\mu}{L + \nicefrac{\cL}{n}}K\right) + \tfrac{\sigma_*^2}{\mu n K} + \tfrac{L\Delta^2}{\mu^2 K^2}\right),&& \text{when } \mu > 0, \label{eq:complexity_str_cvx_EC_SGD_AS}\\
        \cO\left(\tfrac{(L + \nicefrac{\cL}{n})R_0^2}{K} + \sqrt{\tfrac{\sigma_*^2R_0^2}{nK}} + \tfrac{L^{\nicefrac{1}{3}}R_0^{\nicefrac{4}{3}}\Delta^{\nicefrac{2}{3}}}{K^{\nicefrac{2}{3}}}\right),&& \text{when } \mu = 0. \label{eq:complexity_cvx_EC_SGD_AS}
    \end{eqnarray}
\end{theorem}

Consider the case when $\mu > 0$ (similar observations are valid when $\mu = 0$). In these settings, Assumption~\ref{as:expected_smoothness} is satisfied whenever $f_{ij}$ are convex and smooth without assuming $(M,\sigma^2)$-bounded noise assumption used in \cite{sahu2021rethinking}. Moreover, our bound \eqref{eq:complexity_str_cvx_EC_SGD_AS} has better $\cO(\nicefrac{1}{K})$ decaying term than bound \eqref{eq:complexity_str_cvx_EC_SGD_M_sigma} from \cite{sahu2021rethinking}. In particular, when $\sigma_*^2 = 0$, i.e., workers compute full gradients our bound has $\cO(\nicefrac{1}{K^2})$ decaying leading term while for \eqref{eq:complexity_str_cvx_EC_SGD_M_sigma} the leading term decreases as $\cO(\nicefrac{1}{K})$, when $\zeta_*^2 = 0$ (local functions has different optima). Next, when $\mu > 0$, the best-known bound for \algname{EC-SGD-AS} for $\delta$-contractive compressors is (see \cite{gorbunov2020linearly})
\begin{eqnarray}
    \widetilde{\cO}\left(AR_0^2\exp\left(-\tfrac{\mu}{A}K\right) + \tfrac{\sigma_*^2}{\mu n K} + \tfrac{L(\sigma_*^2 + \nicefrac{\zeta_*^2}{\delta})}{\delta \mu^2 K^2}\right), \label{eq:complexity_str_cvx_EC_SGD_AS_contractive}
\end{eqnarray}
where $A = \cL + \tfrac{\max_{i\in [n]}L_i + \sqrt{\delta \max_{i\in[n]}L_i\cL}}{\delta}$ (this parameter can be tightened using the independence of the samples on different workers, which we use in our proofs). The second terms from \eqref{eq:complexity_str_cvx_EC_SGD_AS} and \eqref{eq:complexity_str_cvx_EC_SGD_AS_contractive} are the same while the third terms are different. Although, these results are derived for different classes of compressors, one can compare them for particular choices of compressions. In particular, for hard-threshold and Top1 compressors the third term in \eqref{eq:complexity_str_cvx_EC_SGD_AS} is proportional to $\nicefrac{d \lambda}{K^2}$, while the corresponding term from \eqref{eq:complexity_str_cvx_EC_SGD_AS_contractive} is proportional to $\nicefrac{(d\sigma_*^2 + d^2\zeta_*^2)}{K^2}$. When $\lambda = \cO(1)$ (e.g., see Fig.~\ref{fig:exp_results}) and $\zeta_*^2$ is large enough the bound \eqref{eq:complexity_str_cvx_EC_SGD_AS} is more than $d$ times better than \eqref{eq:complexity_str_cvx_EC_SGD_AS_contractive}.

\section{Absolute Compression and Variance Reduction}\label{sec:lsvrg}

In the same setup as in the previous section, we consider a variance-reduced version of \algname{EC-SGD} called Error Compensated Loopless Stochastic Variance-Reduced Gradient (\algname{EC-LSVRG}) from \cite{gorbunov2020linearly}. This method is a combination of \algname{LSVRG} \cite{hofmann2015variance,kovalev2019don} and \algname{EC-SGD} and can be viewed as Algorithm~\ref{alg:EC_SGD} with
\begin{eqnarray}
    g_i^k &=& \nabla f_{\xi_i^k}(x^k) - \nabla f_{\xi_i^k}(w^k) + \nabla f_i(w^k) \label{eq:LSVRG_g_i^k}\\
    w^{k+1} &=& \begin{cases} x^k, & \text{with probability } p,\\ w^k, & \text{with probability } 1-p, \end{cases}\quad w^0 = x^0, \label{eq:LSVRG_w_i^k}
\end{eqnarray}
where $\xi_i^k$ is sampled from $\cD_i$ independently from previous steps and other workers, and probability $p$ of updating $w^k$ is usually taken as $p \sim \nicefrac{1}{m}$. Such choice of $p$ ensures that full gradients $\nabla f_i(w^k)$ are computed rarely meaning that the expected number of $\nabla f_{\xi_i}(x)$ computations per iteration is the same as for \algname{EC-SGD-AS} (up to the constant factor). We point out that \algname{EC-LSVRG} was studied for the contractive compressors and uniform sampling \cite{gorbunov2020linearly} (although it is possible to generalize the proofs from \cite{gorbunov2020linearly} to cover \algname{EC-LSVRG} with arbitrary sampling as well). As we show next, \algname{EC-LSVRG} with arbitrary sampling satisfies Assumption~\ref{as:key_assumption} (the proof is deferred to Appendix~\ref{appendix:proofs_special_cases}).
\begin{proposition}\label{prop:EC_LSVRG_fits_assumption}
    Let $f$ be $L$-smooth, $f_i$ have finite-sum form, and Assumption~\ref{as:expected_smoothness} hold. Then the iterates produced by \algname{EC-LSVRG} satisfy Assumption~\ref{as:key_assumption} with $A = L + \tfrac{2\cL}{n}$, $B = \tfrac{2}{n}$, $D_1 = 0$, $\sigma_k^2 = 2\cL(f(w^k) - f(x^*))$, $\rho = p$, $C = p\cL$, $D_2 = 0$.
\end{proposition}

Due to the variance reduction, noise terms $D_1$ and $D_2$ equal zero for \algname{EC-LSVRG} allowing the method to achieve better accuracy with constant stepsize than \algname{EC-SGD-AS}. Plugging the parameters from the above proposition in Theorem~\ref{thm:main_result} and Corollary~\ref{cor:main_result}, one can derive convergence guarantees for \algname{EC-LSVRG} with absolute compression operator. In particular, our general analysis implies the following result.
\begin{theorem}\label{thm:EC_LSVRG_convergence}
   Let the assumptions of Proposition~\ref{prop:EC_LSVRG_fits_assumption} hold. Then, there exist choices of stepsize $0 < \gamma \leq \gamma_0 = (4L + \nicefrac{152\cL}{3n})^{-1}$ such that \algname{EC-LSVRG} with absolute compressor guarantees $\EE[f(\overline{x}^K) - f(x^*)]$ of the order\footnote{We take into account that due to $L$-smoothness of $f$ we have $T_0 = \|x^0 - x^*\|^2 + \tfrac{64m}{3n}\gamma^2\cL(f(x^0) - f(x^*)) \leq (1 + \tfrac{32m\cL L\gamma^2}3n)\|x^0 - x^*\|^2$.}
    \begin{eqnarray}
        \widetilde{\cO}\left(\!\left(L + \tfrac{\cL}{n}\right)\widetilde{T}_0\exp\left(\!-\min\left\{\tfrac{\mu}{L + \nicefrac{\cL}{n}}, \tfrac{1}{m}\right\}K\!\right) + \tfrac{L\Delta^2}{\mu^2 K^2}\!\right)\!,&& \text{when } \mu > 0, \label{eq:complexity_str_cvx_EC_LSVRG}\\
        \cO\left(\!\tfrac{(L + \nicefrac{\cL}{n})R_0^2}{K} + \tfrac{\sqrt{m\cL L}R_0^2}{\sqrt{n}K} + \tfrac{L^{\nicefrac{1}{3}}R_0^{\nicefrac{4}{3}}\Delta^{\nicefrac{2}{3}}}{K^{\nicefrac{2}{3}}}\!\right)\!,&& \text{when } \mu = 0, \label{eq:complexity_cvx_EC_LSVRG}
    \end{eqnarray}
    where $\widetilde{T}_0 = \|x^0 - x^*\|^2 + \tfrac{64m}{3n}\gamma_0^2\cL(f(x^0) - f(x^*))$.
\end{theorem}

Consider the case when $\mu > 0$ (similar observations are valid when $\mu = 0$). As expected, the bound \eqref{eq:complexity_str_cvx_EC_LSVRG} does not have terms proportional to any kind of variance of the stochastic estimator. Therefore, the leading term in the complexity bound for \algname{EC-LSVRG} decreases as $\cO(\nicefrac{1}{K^2})$, while \algname{EC-SGD-AS} has $\cO(\nicefrac{1}{K})$ leading term. Next, in case of $\delta$-contractive compression operator the only known convergence rate for \algname{EC-LSVRG} \cite{gorbunov2020linearly} has the leading term $\cO(\tfrac{L\zeta_*^2}{\delta^2\mu^2 K^2})$. Using the same arguments as in the discussion after Theorem~\ref{thm:EC_SGD_AS_convergence}, one can that the leading term in the case of hard-threshold sparsifier can be more than $d$ times better than the leading term in the case of Top1 compressor.

\section{Numerical Experiments}\label{sec:num_exp}
\begin{figure}[t]
    \centering
    \includegraphics[width=0.325\textwidth]{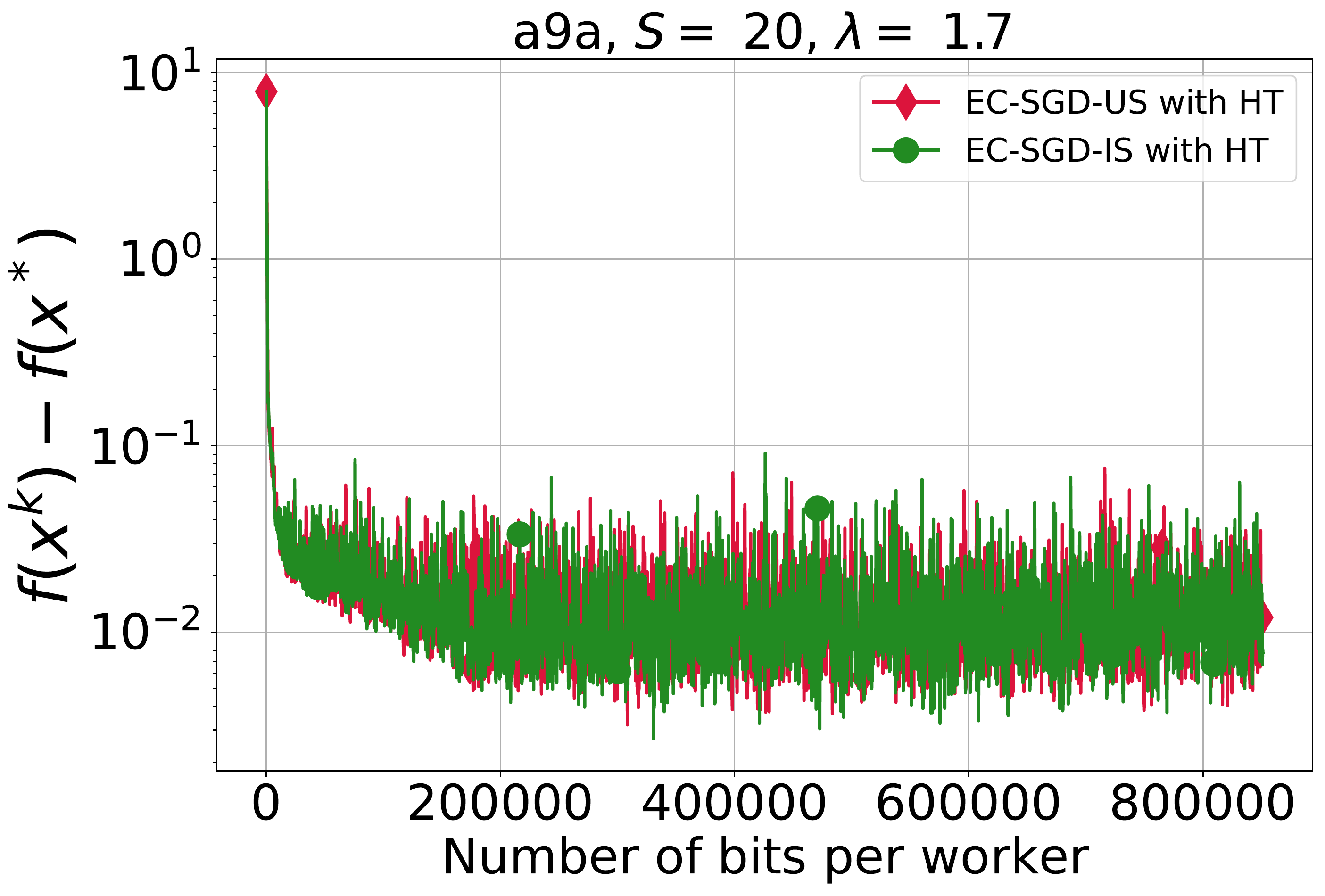}
    \includegraphics[width=0.325\textwidth]{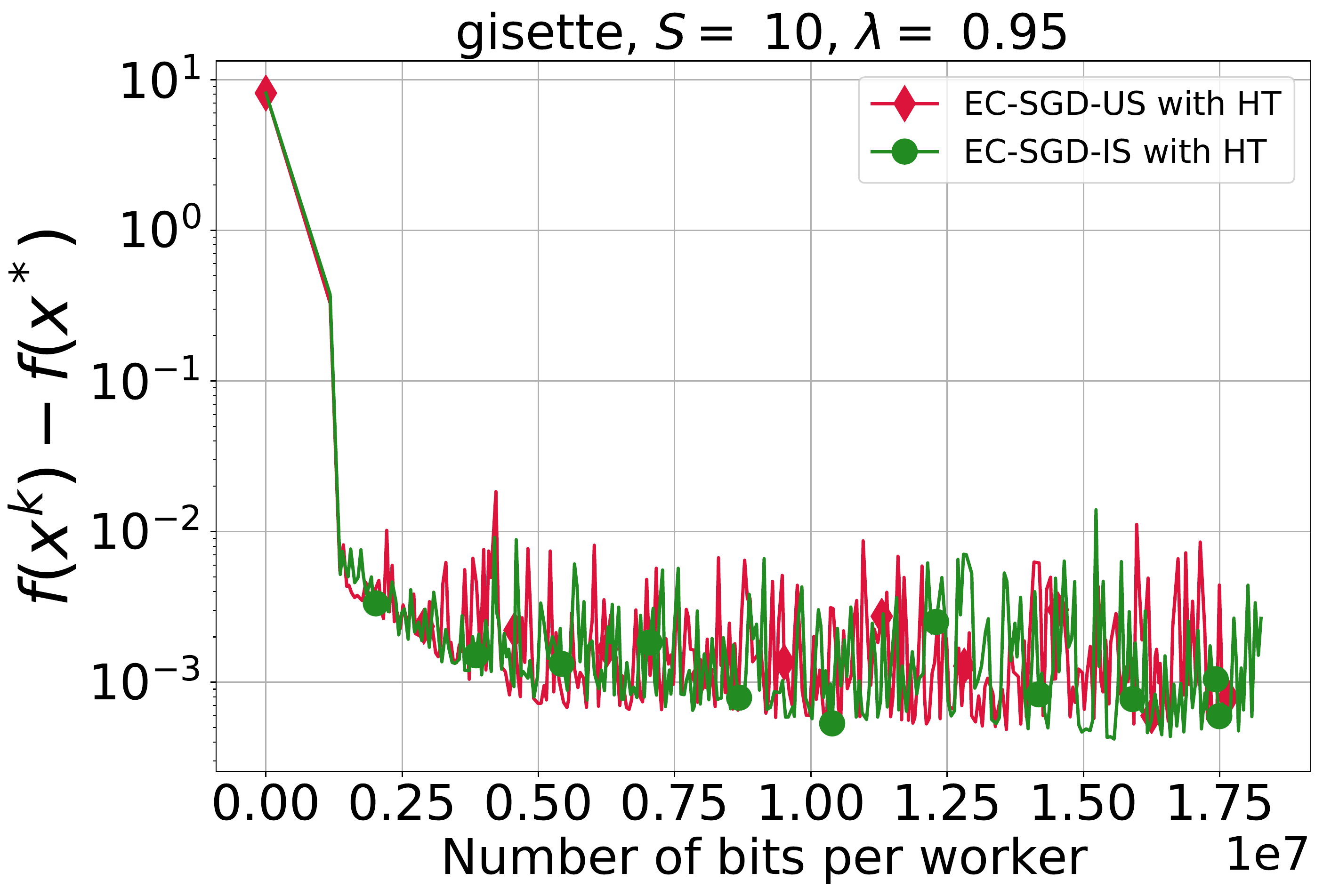}
    \includegraphics[width=0.325\textwidth]{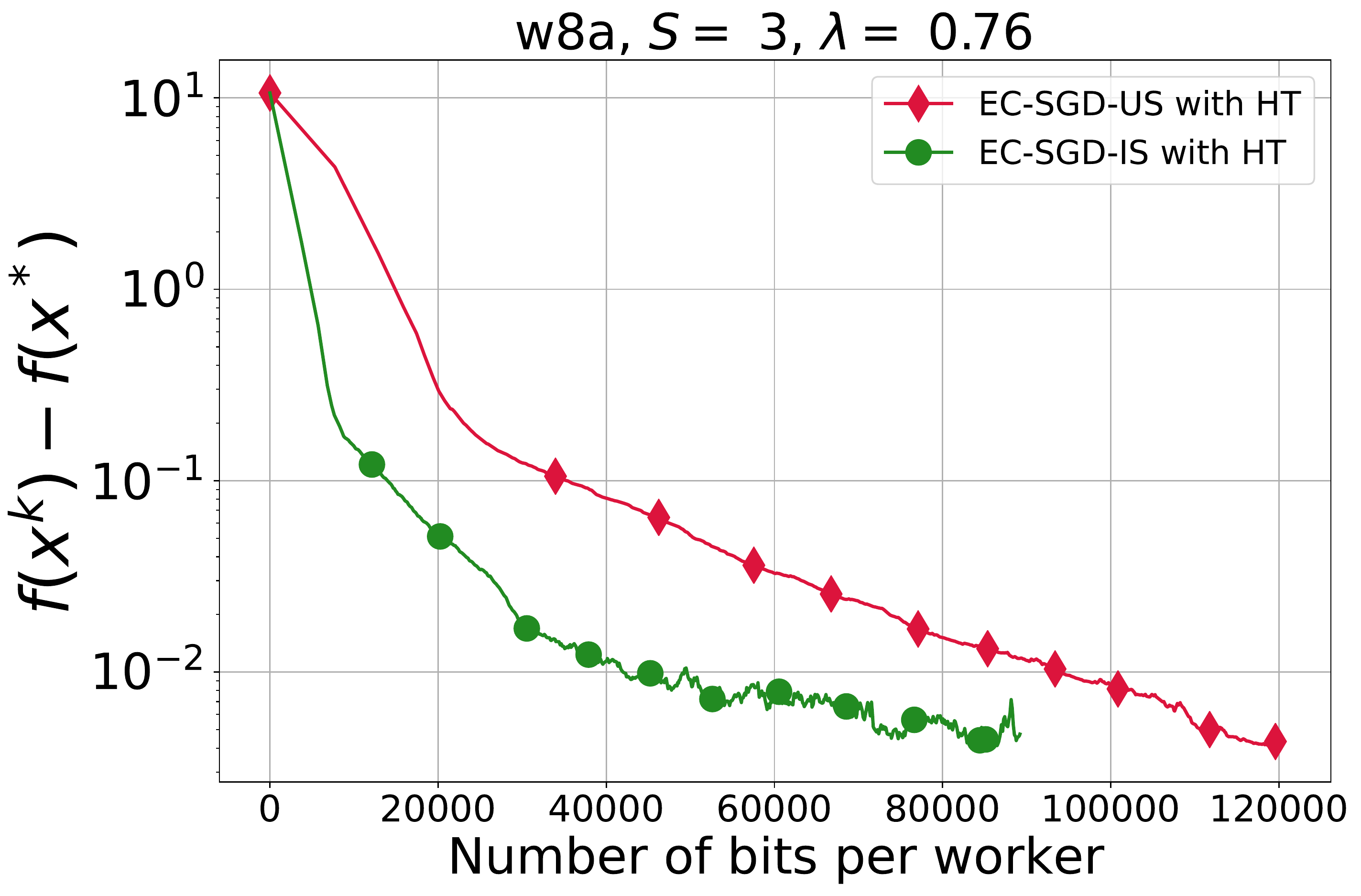}
    \includegraphics[width=0.325\textwidth]{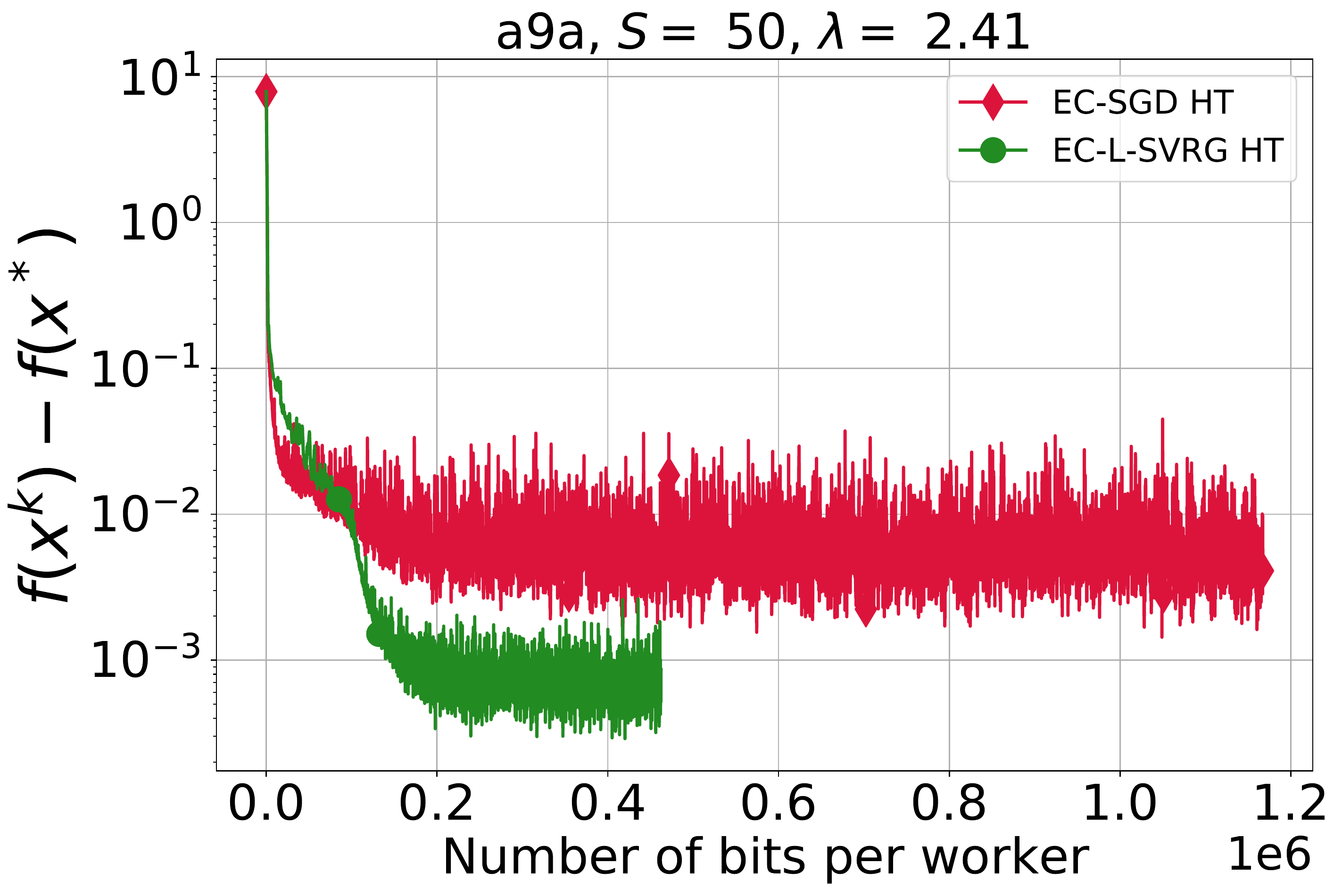}
    \includegraphics[width=0.325\textwidth]{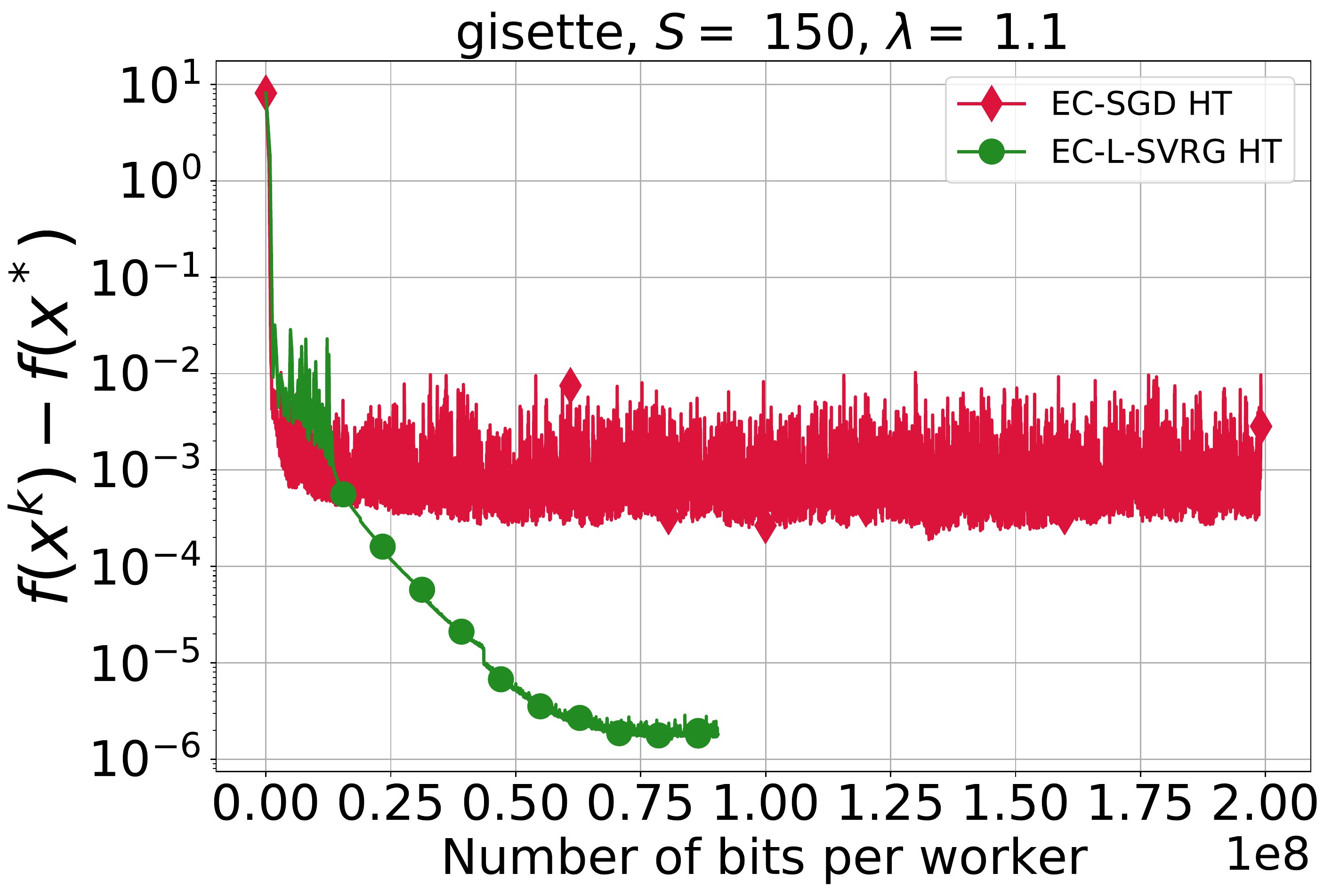}
    \includegraphics[width=0.325\textwidth]{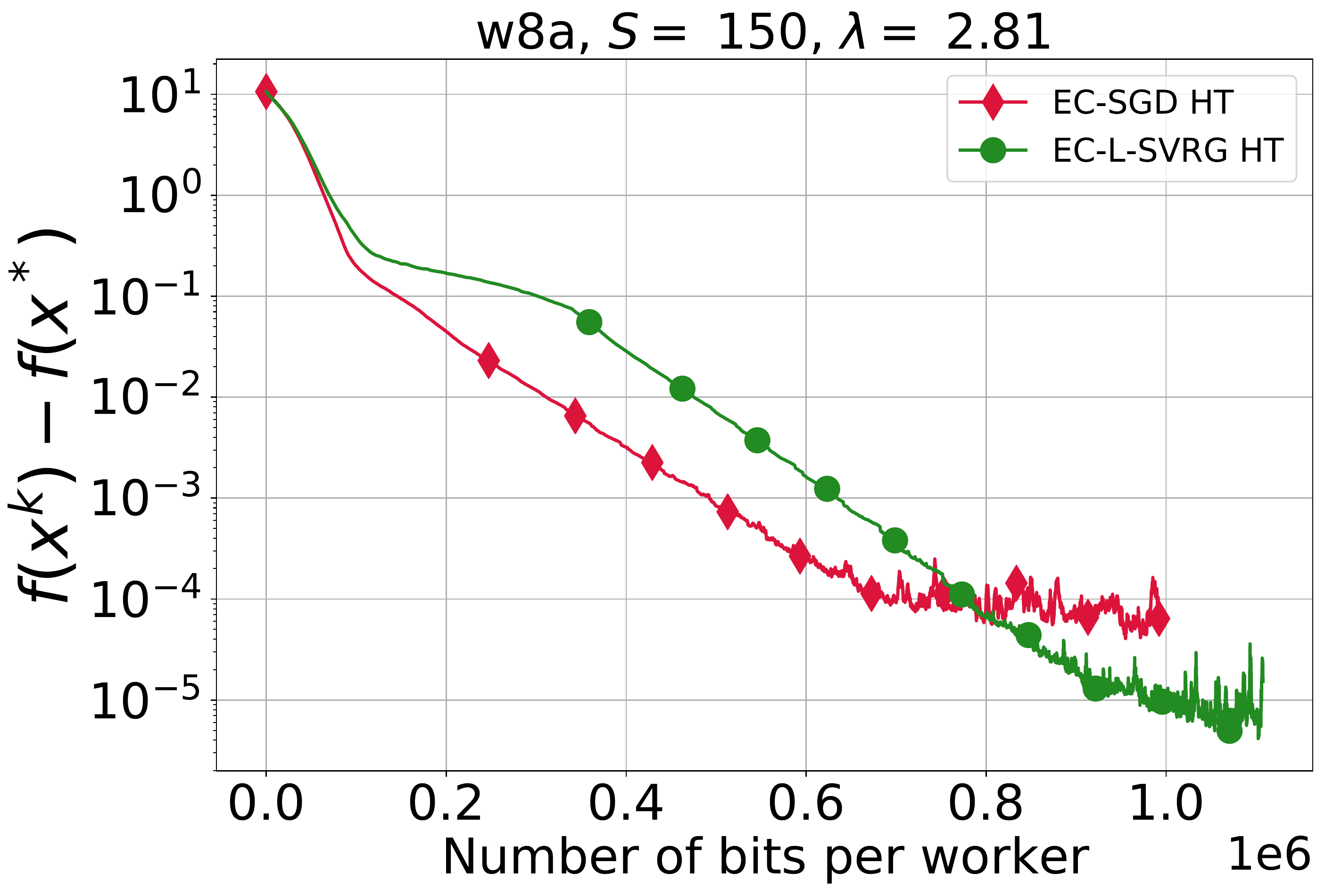}
    \includegraphics[width=0.325\textwidth]{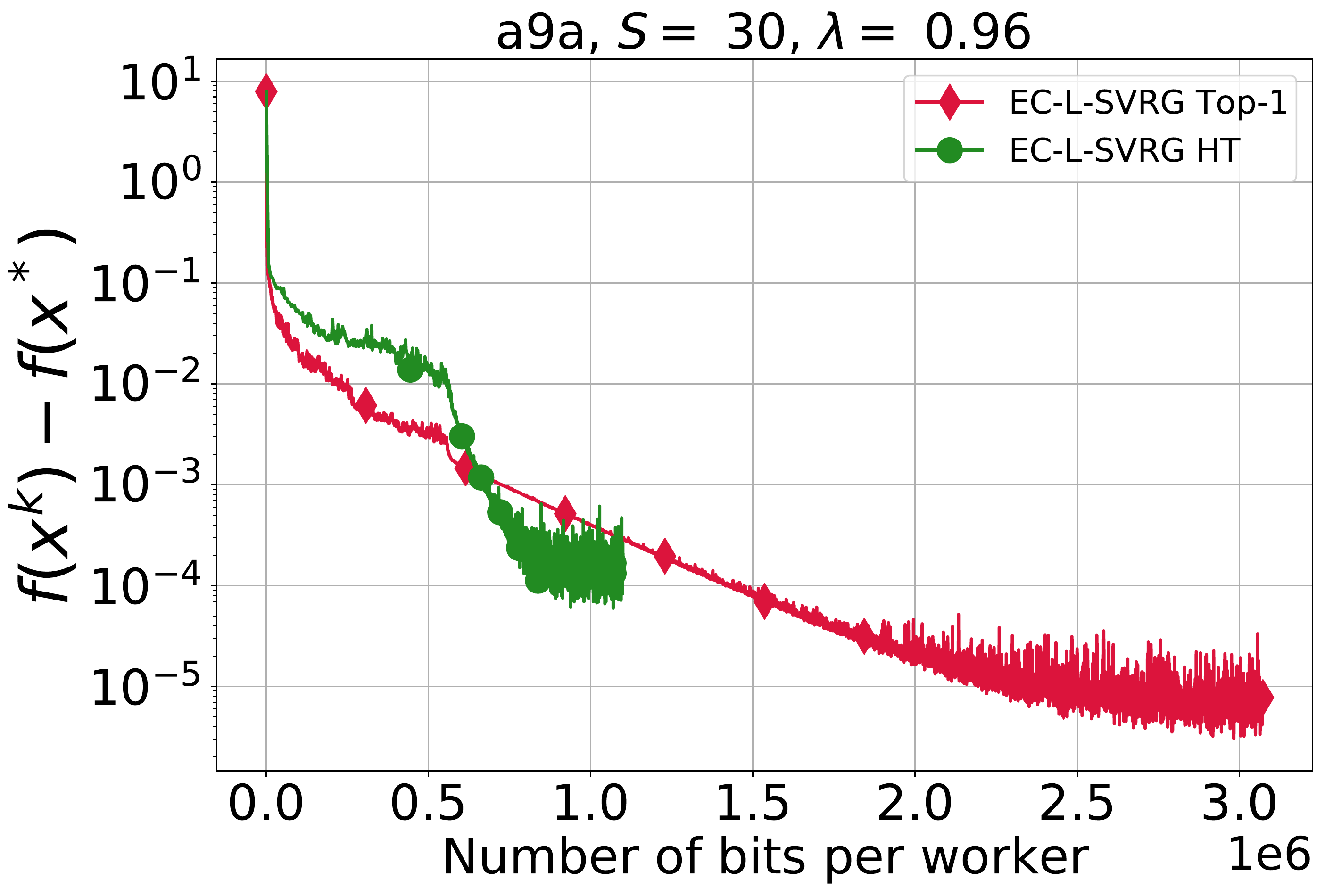}
    \includegraphics[width=0.325\textwidth]{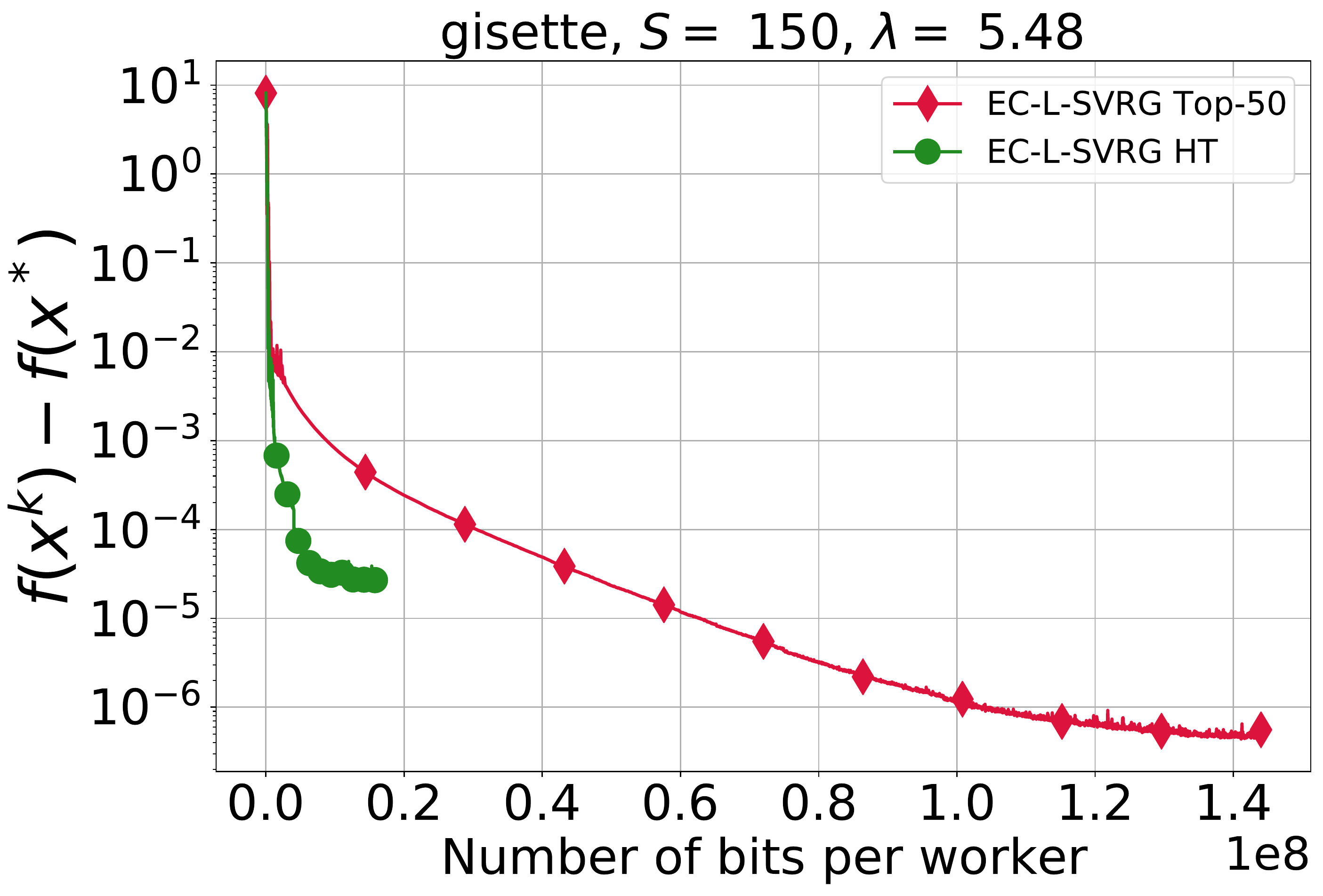}
    \includegraphics[width=0.325\textwidth]{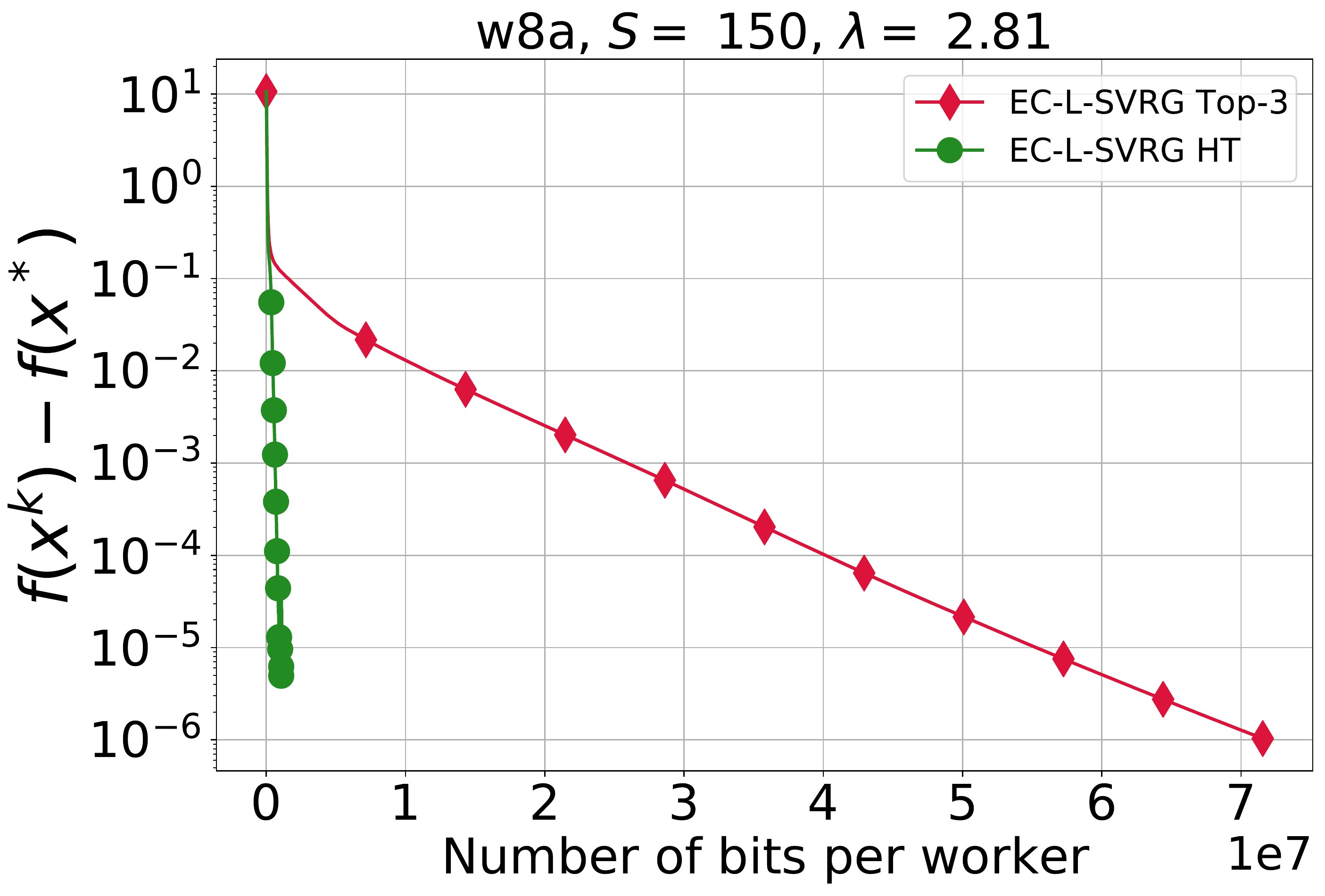}
    \caption{Trajectories of \algname{EC-SGD} with uniform and importance samplings sampling and hard-threshold sparsifier, \algname{EC-LSVRG} with hard-threshold and TopK sparsifiers. The row $i$ of plots corresponds to Experiment $i$ described in Section~\ref{sec:num_exp}.}
    \label{fig:exp_results}
\end{figure}

We conduct several numerical experiments to support our theory, i.e., we tested the methods on the distributed logistic regression problem with $\ell_2$-regularization:
\begin{equation}
    \min\limits_{x\in \R^d}\Bigg\{f(x) = \tfrac{1}{nm}\sum_{i=1}^n\sum_{j=1}^m \underbrace{\ln\left(1 + \exp\left(-y_i \langle a_{ij}, x\rangle\right)\right) + \tfrac{l_2}{2}\|x\|^2}_{f_{ij}(x)}\Bigg\},  \label{eq:logreg}
\end{equation}
where vectors $\{a_{ij}\}_{i\in [n], j\in [m]} \in \R^d$ are the columns of the matrix of features $\mA^\top$, $\{y_{i}\}_{i=1}^n \in \{-1,1\}$ are labels, and $l_2 \geq 0$ is a regularization parameter. One can show that $f_{ij}$ is $l_2$-strongly convex and $L_{ij}$-smooth, and $f$ is $L$-smooth with $L_{ij} = l_2 + \nicefrac{\|a_{ij}\|^2}{4}$ and $L = l_2 + \nicefrac{\lambda_{\max}(\mA^\top\mA)}{4nm}$, where $\lambda_{\max}(\mA^\top\mA)$ is the largest eigenvalue of $\mA^\top\mA$. In particular, we took 3 datasets -- {\tt a9a} ($nm = 32000$, $d = 123$)\footnote{We take the first $32000$ and $49700$ samples from {\tt a9a} and {\tt w8a} to get a multiple of $n = 20$.}, {\tt gisette} ($nm = 6000$, $d = 5000$), and {\tt w8a} ($nm = 49700$, $d = 300$) -- from LIBSVM library \cite{chang2011libsvm}. Out code is based on the one from \cite{gorbunov2020linearly}. Each dataset was shuffled and equally split between $n = 20$ workers. In the first two experiments, we use only hard-threshold sparsifier and in the third experiment, we also consider TopK compressor. The results are presented in Fig.~\ref{fig:exp_results}. The methods were run for $S$ epochs, where the values of $S$ are given in the titles of corresponding plots. We compare the methods in terms of total number of bits each worker sends to the server (on average).

\paragraph{Experiment 1: \algname{EC-SGD} with and without importance sampling.} In this experiment, we tested \algname{EC-SGD} with hard-threshold sparsifier and two different sampling strategies: uniform sampling (US) and importance sampling (IS), described in Section~\ref{sec:arb_sampling}. We chose $l_2 = 10^{-4}\cdot\max_{i\in [n]}\overline{L}_i$. Stepsize was chosen as $\gamma_{\text{US}} = (L + n^{-1}\max_{i\in[n], j\in[m]}L_{ij})^{-1}$ and $\gamma_{\text{IS}} = (L + n^{-1}\max_{i\in[n]}\overline{L}_{i})^{-1}$ for the case of US and IS respectively, which are the multiples of the maximal stepsizes that our theory allows for both cases. Following \cite{sahu2021rethinking}, we took parameter $\lambda$ as $\lambda = 5000\sqrt{\nicefrac{\varepsilon}{d^2 \gamma_{US}}}$ for $\varepsilon = 10^{-3}$. We observe that \algname{EC-SGD} behaves similarly in both cases for {\tt a9a} ($L \approx 1.57$, $\max_{i\in [n]}\overline{L}_i \approx  3.47$, $\max_{i\in[n], j\in[m]}L_{ij} \approx 3.5$) and {\tt gisette} ($L \approx 842.87$, $\max_{i\in [n]}\overline{L}_i \approx  1164.89$, $\max_{i\in[n], j\in[m]}L_{ij} \approx 1201.51$), while for {\tt w8a} ($L \approx 0.66$, $\max_{i\in [n]}\overline{L}_i \approx  3.05$, $\max_{i\in[n], j\in[m]}L_{ij} \approx 28.5$) \algname{EC-SGD} with IS achieves good enough accuracy much faster than with US. This is expected since for {\tt w8a} $\max_{i\in [n]}\overline{L}_i$ is almost $10$ times smaller than $\max_{i\in[n], j\in[m]}L_{ij}$. That is, as our theory implies, importance sampling is preferable when $\max_{i\in [n]}\overline{L}_i \ll \max_{i\in[n], j\in[m]}L_{ij}$.

\paragraph{Experiment 2: \algname{EC-SGD} vs \algname{EC-LSVRG}.} Next, we compare \algname{EC-SGD} and \algname{EC-LSVRG} to illustrate the benefits of variance reduction for error-compensated methods with absolute compression. Both methods were run with stepsize $\gamma = \nicefrac{1}{\max_{i\in[n], j\in[m]}L_{ij}}$ and $\lambda = 5000\sqrt{\nicefrac{\varepsilon}{d^2 \gamma}}$ for $\varepsilon = 10^{-3}$. In all cases, \algname{EC-LSVRG} achieves better accuracy of the solution than \algname{EC-SGD} that perfectly corroborates our theoretical results.

\paragraph{Experiment 3: \algname{EC-LSVRG} with hard-threshold and TopK sparsifiers.} Finally, to highlight the benefits of hard-threshold (HT) sparsifier compared to TopK sparsifier we tested \algname{EC-LSVRG} with both compressors. We chose stepsize as $\gamma = \nicefrac{1}{\max_{i\in[n], j\in[m]}L_{ij}}$, $K \approx \nicefrac{d}{100}$ (see the legends of plots in row 3 from Fig.~\ref{fig:exp_results}), and $\lambda = \alpha\sqrt{\nicefrac{\varepsilon}{d^2 \gamma}}$, where $\varepsilon = 10^{-3}$ and $\alpha = 2000$ for {\tt a9a}, $\alpha = 25000$ for {\tt gisette}, $\alpha = 5000$ for {\tt w8a}. We observe that \algname{EC-LSVRG} with HT achieves a reasonable accuracy ($10^{-3}$ for {\tt a9a}, $10^{-4}$ for {\tt gisette}, and $10^{-5}$ for {\tt w8a}) faster than \algname{EC-LSVRG} with TopK for all datasets. In particular,  \algname{EC-LSVRG} with HT significantly outperforms \algname{EC-LSVRG} with TopK on {\tt w8a} dataset.

\bibliography{refs}

\appendix

\section{Missing Proofs from Section~\ref{sec:unif_analysis}}\label{appendix:main_proofs}
In the analysis, we use auxiliary iterates that are never computed explicitly during the work of the method: $\tx^{k} = x^k - e^k$, where $e^k = \tfrac{1}{n}\sum_{i=1}^n e_i^k$. 
These iterates are usually called \emph{perturbed} or \emph{virtual} iterates \cite{leblond2018improved,mania2017perturbed}. They are used in many previous works on error feedback \cite{stich2018sparsified,stich2020error,gorbunov2020linearly,sahu2021rethinking}. One of the key properties these iterates satisfy is the following recursion:
\begin{equation}
    \tx^{k+1} = x^{k+1} - e^{k+1} = x^k - v^k - e^k - \gamma g^k + v^k = \tx^k - \gamma g^k,\label{eq:virtual_it_recursion}
\end{equation}
where we use $e^{k+1} = e^k +\gamma g^k - v^k$, which follows from $e_i^{k+1} = e_i^k +\gamma g_i^k - v_i^k$ and definitions of $e^k, g^k$, and $v^k$.

\subsection{Proof of Theorem~\ref{thm:main_result}}
Our proof is close to the ones from \cite{gorbunov2020linearly,sahu2021rethinking}. Using the recursion \eqref{eq:virtual_it_recursion} for the perturbed iterates, we obtain
\begin{eqnarray*}
    \|\tx^{k+1} - x^*\|^2 &=& \|\tx^k - x^*\|^2 - 2\gamma \langle \tx^k - x^*, g^k \rangle + \gamma^2 \|g^k\|^2.
\end{eqnarray*}
Next, we take conditional expectation $\EE_k[\cdot]$ from the above inequality and apply unbiasedness of $g^k$ and inequality \eqref{eq:second_moment_bound_general} from Assumption~\ref{as:key_assumption}:
\begin{eqnarray*}
    \EE_k\left[\|\tx^{k+1} - x^*\|^2\right] 
    &\overset{\eqref{eq:second_moment_bound_general}}{\leq}& 
    \|\tx^k - x^*\|^2 - 2\gamma \langle x^k - x^*, \nabla f(x^k) \rangle + 2\gamma \langle x^k - \tx^k, \nabla f(x^k) \rangle\\
    &&\quad + \gamma^2\left(2A\left(f(x^k) - f(x^*)\right) + B\sigma_k^2 + D_1\right).
\end{eqnarray*}
Since $f$ is $\mu$-quasi strongly convex, we have $-\langle x^k - x^*, \nabla f(x^k) \rangle \overset{\eqref{eq:quasi_strong_cvx}}{\leq} -(f(x^k) - f(x^*)) - \tfrac{\mu}{2}\|x^k - x^*\|^2$ implying
\begin{eqnarray}
    \EE_k\left[\|\tx^{k+1} - x^*\|^2\right] &=& \|\tx^k - x^*\|^2 - \gamma\mu \|x^k - x^*\|^2 + 2\gamma \langle x^k - \tx^k, \nabla f(x^k) \rangle \notag\\
    &&\quad -2\gamma(1 - A\gamma)\left(f(x^k) - f(x^*)\right) + B\gamma^2\sigma_k^2 + \gamma^2D_1 \label{eq:main_proof_technical_1}
\end{eqnarray}
To continue our derivation we need to upper bound $- \gamma\mu \|x^k - x^*\|^2$ and $2\gamma \langle x^k - \tx^k, \nabla f(x^k) \rangle$. Applying $\|a-b\|^2 \geq \tfrac{1}{2}\|a\|^2 - \|b\|^2$, which holds for any $a,b \in \R^d$, we get
\begin{eqnarray}
    - \gamma\mu \|x^k - x^*\|^2 
    &\leq& -\tfrac{\gamma\mu}{2}\|\tx^k - x^*\|^2 + \gamma\mu\|\tx^k - x^k\|^2. \label{eq:main_proof_technical_2}
\end{eqnarray}
To estimate the inner product, we use Fenchel-Young inequality $\langle a, b\rangle \leq \tfrac{\alpha}{2}\|a\|^2 + \tfrac{1}{2\alpha}\|b\|^2$ holding for any $a,b\in\R^d$, $\alpha > 0$ together with standard inequality $\|\nabla f(x^k)\|^2 \leq 2L(f(x^k) - f(x^*))$, which holds for any $L$-smooth function $f$ \cite{nesterov2018lectures}:
\begin{eqnarray}
    2\gamma \langle x^k - \tx^k, \nabla f(x^k) \rangle &\leq& 2\gamma L \|x^k - \tx^k\|^2 + \tfrac{\gamma}{2L}\|\nabla f(x^k)\|^2 \notag\\
    &\leq& 2\gamma L \|\tx^k - x^k\|^2 + \gamma\left(f(x^k) - f(x^*)\right). \label{eq:main_proof_technical_3}
\end{eqnarray}
Plugging upper bounds \eqref{eq:main_proof_technical_2} and \eqref{eq:main_proof_technical_3} in \eqref{eq:main_proof_technical_1}, we derive
\begin{eqnarray}
    \EE_k\left[\|\tx^{k+1} - x^*\|^2\right] &\leq& \left(1 - \tfrac{\gamma\mu}{2}\right)\|\tx^k - x^*\|^2 - \gamma(1 - 2A\gamma)\left(f(x^k) - f(x^*)\right) \notag\\
    &&\quad + B\gamma^2\sigma_k^2 + \gamma^2D_1  + \gamma(2L + \mu) \|x^k - \tx^k\|^2 \notag\\
    &\overset{\mu \leq L}{\leq}&  \left(1 - \tfrac{\gamma\mu}{2}\right)\|\tx^k - x^*\|^2 - \gamma(1 - 2A\gamma)\left(f(x^k) - f(x^*)\right) \notag\\
    &&\quad + B\gamma^2\sigma_k^2 + \gamma^2D_1  + 3L\gamma \|e^k\|^2. \notag
\end{eqnarray}
Since the definition of $e^{k+1}$ and Jensen's inequality imply that for all $k\ge 0$
\begin{eqnarray*}
    \|e^{k+1}\|^2 &\leq& \tfrac{1}{n}\sum_{i=1}^n \|e_i^k + \gamma g_i^k - v_i^k\|^2 = \tfrac{\gamma^2}{n}\sum_{i=1}^n \left\|\tfrac{e_i^k + \gamma g_i^k}{\gamma} - \cC\left(\tfrac{e_i^k + \gamma g_i^k}{\gamma}\right)\right\|^2 \leq \gamma^2\Delta^2,
\end{eqnarray*}
we have (taking into account that $\|e^0\|^2 = 0 \leq \gamma^2\Delta^2$)
\begin{eqnarray}
    \EE_k\left[\|\tx^{k+1} - x^*\|^2\right] &\leq& \left(1 - \tfrac{\gamma\mu}{2}\right)\|\tx^k - x^*\|^2 - \gamma(1 - 2A\gamma)\left(f(x^k) - f(x^*)\right) \notag\\
    &&\quad + B\gamma^2\sigma_k^2 + \gamma^2D_1  + 3L\gamma^3 \Delta^2. \notag
\end{eqnarray}
Summing up the above inequality with $F\gamma^2$-multiple of \eqref{eq:sigma_k_bound_geenral} and introducing new notation $T_k = \|\tx^k - x^*\|^2 + F\gamma^2 \sigma_k^2$, we obtain
\begin{eqnarray*}
    \EE_k\left[T_{k+1}\right] &\leq& \left(1 - \tfrac{\gamma\mu}{2}\right)\|\tx^k - x^*\|^2 + F\gamma^2\left(1 - \rho + \tfrac{B}{F}\right)\sigma_k^2\\
    &&\quad - \gamma(1 - 2(A+CF)\gamma)\left(f(x^k) - f(x^*)\right)\\
    &&\quad + \gamma^2 \left(D_1 + FD_2 + 3L\gamma \Delta^2\right)\\
    &\overset{F = \nicefrac{4B}{3\rho}}{\leq}& (1-\eta) T_k + \gamma^2 \left(D_1 + FD_2 + 3L\gamma \Delta^2\right) - \tfrac{\gamma}{2}\left(f(x^k) - f(x^*)\right),
\end{eqnarray*}
where in the last step we use $\left(1 - \tfrac{\gamma\mu}{2}\right)\|\tx^k - x^*\|^2 + \left(1 - \tfrac{\rho}{4}\right)F\gamma^2\sigma_k^2 \leq (1 - \eta)\left(\|\tx^k - x^*\|^2 + F\gamma^2 \sigma_k^2\right) = (1-\eta)T_k$, where $\eta = \min\{\nicefrac{\gamma\mu}{2}, \nicefrac{\rho}{4}\}$, and $0 < \gamma \leq \nicefrac{1}{4(A+CF)}$. Rearranging the terms and taking full expectation, we derive
\begin{equation*}
    \tfrac{\gamma}{2}\EE\left[f(x^k) - f(x^*)\right] \leq (1-\eta) \EE[T_k] - \EE[T_{k+1}] + \gamma^2 \left(D_1 + FD_2 + 3L\gamma \Delta^2\right).
\end{equation*}
Summing up the above inequality for $k = 0,\ldots, K$ with weights $w_k = (1-\eta)^{-(k+1)}$ and using Jensen's inequality $f(\overline{x}^K) \leq \tfrac{1}{W_K}\sum_{k=0}^K w_k f(x^k)$ for convex function $f$ and point $\overline{x}^K = \tfrac{1}{W_K}\sum_{k=0}^K w_k x^k$, where $W_K = \sum_{k=0}^K w_k$, and $w_{k} = (1-\eta)w_{k-1}$, we get
\begin{eqnarray}
    \EE\left[f(\overline{x}^K) - f(x^*)\right] &\leq& \tfrac{2}{\gamma W_K}\sum_{k=0}^K\left(w_{k-1}\EE[T_k] - w_k\EE[T_{k+1}]\right)\notag \\
    &&\quad + 2\gamma \left(D_1 + FD_2 + 3L\gamma \Delta^2\right)\notag \\
    &\leq& \tfrac{2\EE[T_0]}{\gamma W_k} + 2\gamma \left(D_1 + FD_2 + 3L\gamma \Delta^2\right), \label{eq:main_proof_technical_4}
\end{eqnarray}
where in the last step we use $\sum_{k=0}^K\left(w_{k-1}\EE[T_k] - w_k\EE[T_{k+1}]\right) = w_{-1}\EE[T_0] - w_{k+1}\EE[T_{K+1}] \leq w_{-1}\EE[T_0] = \EE[T_0]$. Finally, it remains to notice that \eqref{eq:main_proof_technical_4} implies \eqref{eq:main_result_str_cvx} and \eqref{eq:main_result_cvx}. Indeed, when $\mu > 0$, we have $W_K \geq w_K = (1 - \eta)^{-(K+1)}$, and when $\mu = 0$, we have $W_K = K+1$, since $\eta = 0$.

\section{Missing Proofs from Sections~\ref{sec:arb_sampling} and \ref{sec:lsvrg}}\label{appendix:proofs_special_cases}

\subsection{Proof of Proposition~\ref{prop:EC_SGD_AS_fits_assumption}}
To prove Proposition~\ref{prop:EC_SGD_AS_fits_assumption} we need to derive an upper bound for $\EE_k[\|g^k\|^2]$. Independence of $\xi_1^k,\ldots,\xi_n^k$ for fixed history, variance decomposition, and standard inequality $\|\nabla f(x^k)\|^2 \leq 2L(f(x) - f(x^*))$ \cite{nesterov2018lectures} imply
\begin{eqnarray*}
    \EE_k\left[\|g^k\|^2\right] &=& \EE_k\left[\|g^k - \nabla f(x^k)\|^2\right] + \|\nabla f(x^k)\|^2\\
    &\leq& \EE_k\left[\left\|\tfrac{1}{n}\sum_{i=1}^n \nabla f_{\xi_i^k}(x^k) - \nabla f_i(x^k)\right\|^2\right] + 2L\left(f(x^k) - f(x^*)\right)\\
    &=& \tfrac{1}{n^2}\sum_{i=1}^n \EE_k\left[\|\nabla f_{\xi_i^k}(x^k) - \nabla f_i(x^k)\|^2\right] + 2L\left(f(x^k) - f(x^*)\right)\\
    &\leq& \tfrac{2}{n^2}\sum_{i=1}^n \EE_k\left[\|\nabla f_{\xi_i^k}(x^k) - \nabla f_{\xi_i^k}(x^*) - (\nabla f_i(x^k) - \nabla f_i(x^*))\|^2\right] \\
    &&\quad + \underbrace{\tfrac{2}{n^2}\sum_{i=1}^n \EE_k\left[\|\nabla f_{\xi_i^k}(x^*) - \nabla f_i(x^*)\|^2\right]}_{\nicefrac{2\sigma_*^2}{n}} + 2L\left(f(x^k) - f(x^*)\right),
\end{eqnarray*}
where in the last step we use that $\|a+b\|^2 \leq 2\|a\|^2 + 2\|b\|^2$ for all $a,b\in \R^d$. Since the variance is upper-bounded by the second moment, we have
\begin{eqnarray*}
    \EE_k\left[\|g^k\|^2\right] &\leq&  \tfrac{2}{n^2}\sum_{i=1}^n \EE_k\left[\|\nabla f_{\xi_i^k}(x^k) - \nabla f_{\xi_i^k}(x^*)\|^2\right] + 2L\left(f(x^k) - f(x^*)\right) + \tfrac{2\sigma_*^2}{n}\\
    &\leq& 2\left(L + \tfrac{2\cL}{n}\right)\left(f(x^k) - f(x^*)\right) + \tfrac{2\sigma_*^2}{n},
\end{eqnarray*}
where the second inequality follows from Assumption~\ref{as:expected_smoothness} and $f(x) = \tfrac{1}{n}\sum_{i=1}^n f_i(x)$. The derived inequality implies that Assumption~\ref{as:key_assumption} holds with $A = L + \nicefrac{2\cL}{n}$, $B = 0$, $D_1 = \tfrac{2\sigma_*^2}{n} = \tfrac{2}{n^2}\sum_{i=1}^n \EE[\|\nabla f_{\xi_i}(x^*) - \nabla f_i(x^*)\|^2]$, $\sigma_k^2 \equiv 0$, $\rho = 1$, $C = 0$, $D_2 = 0$.

\subsection{Proof of Proposition~\ref{prop:EC_LSVRG_fits_assumption}}

We start with deriving an upper-bound for $\EE_k[\|g^k\|^2]$. Similarly to the proof of Proposition~\ref{prop:EC_SGD_AS_fits_assumption}, we use independence of $\xi_1^k,\ldots,\xi_n^k$ for fixed history, variance decomposition, and standard inequality $\|\nabla f(x^k)\|^2 \leq 2L(f(x) - f(x^*))$ \cite{nesterov2018lectures}:
\begin{eqnarray*}
    \EE_k\left[\|g^k\|^2\right] &=& \EE_k\left[\|g^k - \nabla f(x^k)\|^2\right] + \|\nabla f(x^k)\|^2\\
    &\leq& \EE_k\left[\left\|\tfrac{1}{n}\sum_{i=1}^n \nabla f_{\xi_i^k}(x^k) - \nabla f_{\xi_i^k}(w^k) + \nabla f_{i}(w^k)   - \nabla f_i(x^k)\right\|^2\right]\\
    &&\quad + 2L\left(f(x^k) - f(x^*)\right)\\
    &=& \tfrac{1}{n^2}\sum_{i=1}^n\EE_k\left[\left\|\nabla f_{\xi_i^k}(x^k) - \nabla f_{\xi_i^k}(w^k) - (\nabla f_i(x^k) - \nabla f_{i}(w^k))\right\|^2\right]\\
    &&\quad + 2L\left(f(x^k) - f(x^*)\right).\\
\end{eqnarray*}
Since the variance is upper-bounded by the second moment and $\|a+b\|^2 \leq 2\|a\|^2 + 2\|b\|^2$ for all $a,b\in \R^d$, we have
\begin{eqnarray}
    \EE_k\left[\|g^k\|^2\right] &\leq& \tfrac{1}{n^2}\sum_{i=1}^n\EE_k\left[\left\|\nabla f_{\xi_i^k}(x^k) - \nabla f_{\xi_i^k}(w^k)\right\|^2\right] + 2L\left(f(x^k) - f(x^*)\right)\notag\\
    &\leq& \tfrac{2}{n^2}\sum_{i=1}^n\EE_k\left[\left\|\nabla f_{\xi_i^k}(x^k) - \nabla f_{\xi_i^k}(x^*)\right\|^2\right]\notag\\
    &&\quad + \tfrac{2}{n^2}\sum_{i=1}^n\EE_k\left[\left\|\nabla f_{\xi_i^k}(w^k) - \nabla f_{\xi_i^k}(x^*)\right\|^2\right] + 2L\left(f(x^k) - f(x^*)\right)\notag\\
    &\leq& 2\left(L + \tfrac{2\cL}{n}\right)\left(f(x^k) - f(x^*)\right) + \underbrace{4\cL\left(f(w^k) - f(x^*)\right)}_{\nicefrac{2\sigma_k^2}{n}}, \label{eq:main_proof_technical_5}
\end{eqnarray}
where the third inequality follows from Assumption~\ref{as:expected_smoothness} and $f(x) = \tfrac{1}{n}\sum_{i=1}^n f_i(x)$. Using the definitions of $\sigma_{k+1}^2$ and $w^{k+1}$, we derive an upper bound for $\EE_k[\sigma_{k+1}^2]$:
\begin{eqnarray}
    \EE_k[\sigma_{k+1}^2] &=& 2\cL\EE_k\left[f(w^{k+1}) - f(x^*)\right] \notag\\
    &=& (1-p)\underbrace{2\cL\left(f(w^k) - f(x^*)\right)}_{\sigma_k^2} + 2p\cL\left(f(x^k) - f(x^*)\right). \label{eq:main_proof_technical_6}
\end{eqnarray}
Inequalities \eqref{eq:main_proof_technical_5} and \eqref{eq:main_proof_technical_6} imply that Assumption~\ref{as:key_assumption} holds with $A = L + \nicefrac{2\cL}{n}$, $B = \nicefrac{2}{n}$, $D_1 = 0$, $\sigma_k^2 = 2\cL(f(w^k) - f(x^*))$, $\rho = p$, $C = p\cL$, $D_2 = 0$.

\section{Auxiliary Results}
We use the following lemmas to derive final convergence rates.

\begin{lemma}[Lemma I.2 from \cite{gorbunov2021local}]\label{lem:str_cvx_complexity}
    Let sequence $\{r_k\}_{k \ge 0}$ satisfy $r_K \leq (1-\eta)^{K}\tfrac{a}{\gamma} + c_1\gamma + c_2\gamma^2$ for all $K \ge 0$, $\eta = \min\{\nicefrac{\gamma\mu}{2}, \nicefrac{\rho}{4}\}$, where $\mu > 0$, $\rho \in (0,1]$, with some constants $a,c_1,c_2 \ge 0$ and $0 < \gamma \leq \nicefrac{1}{h}$. Then for all $K > 0$ such that
    \begin{eqnarray}
        \text{either}&& \tfrac{\ln\left(\max\left\{2, \min\left\{\nicefrac{a\mu^2K^2}{4c_1}, \nicefrac{a\mu^3 K^3}{8c_2}\right\}\right\}\right)}{K} \leq \rho \notag\\
        \text{or}&& \tfrac{1}{h} \leq \tfrac{\ln\left(\max\left\{2, \min\left\{\nicefrac{a\mu^2K^2}{4c_1}, \nicefrac{a\mu^3 K^3}{8c_2}\right\}\right\}\right)}{\mu K} \notag
    \end{eqnarray}
    and for the stepsize
    \begin{equation}
        \gamma = \min\left\{\tfrac{1}{h}, \tfrac{\ln\left(\max\left\{2, \min\left\{\nicefrac{a\mu^2K^2}{4c_1}, \nicefrac{a\mu^3 K^3}{8c_2}\right\}\right\}\right)}{\mu K}\right\} \label{eq:stepsize_str_cvx}
    \end{equation}
    we have that
    \begin{equation}
        r_K = \widetilde{\cO}\left(ha\exp\left(-\min\left\{\tfrac{\mu}{h}, \rho\right\}K\right) + \tfrac{c_1}{\mu K} + \tfrac{c_2}{\mu^2 K^2}\right).
    \end{equation}
\end{lemma}

\begin{lemma}[Lemma I.3 from \cite{gorbunov2021local}]\label{lem:cvx_complexity}
    Let sequence $\{r_k\}_{k \ge 0}$ satisfy $r_K \leq \tfrac{a}{\gamma K} + \tfrac{b\gamma}{K} + c_1\gamma + c_2\gamma$ for all $K > 0$ with some constants $a,b,c_1,c_2 \ge 0$ and $0 < \gamma \leq \nicefrac{1}{h}$. Then for all $K > 0$ and
    \begin{equation}
        \gamma = \min\left\{\tfrac{1}{h}, \sqrt{\tfrac{a}{b}}, \sqrt{\tfrac{a}{c_1 K}}, \sqrt[3]{\tfrac{a}{c_2 K}}\right\} \label{eq:stepsize_cvx}
    \end{equation}
    we have that
    \begin{equation}
        r_K = \cO\left(\tfrac{ha + \sqrt{ab}}{K} + \sqrt{\tfrac{ac_1}{K}} + \tfrac{\sqrt[3]{a^2c_2}}{K^{\nicefrac{2}{3}}}\right).
    \end{equation}
\end{lemma}

\end{document}